\documentclass[12pt]{amsart}
\usepackage{latexsym,amssymb,amsmath,amsthm}
\usepackage{listings}
\usepackage{float}
\usepackage{graphicx}
\usepackage{multicol}
\theoremstyle{definition}

\newcommand{\R}{\mathbb{R}}

\def\A{\mathcal{A}}
\def\G{\mathcal{G}}
\def\V{\mathcal{V}}
\def\tr{\rm{tr}}
\def\Det{\rm{Det}}

\def\dim{\rm{dim}}
\def\ker{\rm{ker}}

\def\deg{{\rm deg}}
\def\sign{{\rm sign}}

\def\per{\rm{per}}
\def\deg{\rm{deg}}
\def\str{\rm{str}}

\title{The Dirac operator of a graph}
\author{Oliver Knill}
\date{June 9, 2013}
\address{
        Department of Mathematics \\
        Harvard University \\
        Cambridge, MA, 02138
        }

\begin{document}
\maketitle
\begin{abstract}
We discuss some linear algebra related to the Dirac matrix 
$D$ of a finite simple graph $G=(V,E)$. 
\end{abstract} 

\section{Introduction}
These are expanded preparation notes to a talk given on June 5, 2013 at ILAS. 
References are in \cite{elemente11,poincarehopf,cherngaussbonnet,knillcalculus,randomgraph,josellisknill,
brouwergraph,knillmckeansinger,indexformula,indexexpectation,randomgraph,diracannouncement,cauchybinet}.
It is a pleasure to thank the organizers of ILAS and especially {\bf Leslie Hogben} and {\bf Louis Deaett}, 
the organizers of the matrices and graph theory minisymposium for their kind invitation to participate. 

\section{Dirac operator}
Given a finite simple graph $G=(V,E)$, let $\G_k$ be the set of $K_{k+1}$ subgraphs and let $\G=\bigcup_{k=0} \G_k$ 
be the union of these sub-simplices. Elements in $\G$ are also called {\bf cliques} in graph theory.
The {\bf Dirac operator} $D$ is a symmetric $v \times v$ matrix $D$, where $v$ is the cardinality of $\G$. 
The matrix $D$ acts on the vector space $\Omega = \oplus_{k=0} \Omega_k$, where 
$\Omega_k$ is the space of scalar functions on $\G_k$. An {\bf orientation} of a $k$-dimensional 
clique $x \sim K_{k+1}$ is a choice of a permutation of its $(k+1)$ vertices. An {\bf orientation of the graph} 
is a choice of orientations on all elements of $\G$. We do not require the graph to be {\bf orientable}; the later would require to have
an orientation which is compatible on subsets or intersections of cliques and as a triangularization of
the Moebius strip demonstrates, the later does not always exist. The {\bf Dirac matrix} $D:\Omega \to \Omega$ is zero 
everywhere except for $D_{ij}=\pm 1$ if $i \subset j$ or $j \subset i$ and $|{\rm dim}(x)-{\rm dim}(y)|=1$ and where the sign
tells whether the induced orientations of $i$ and $j$ match or not. The nonnegative matrix $|D_{ij}|$,when  seen as an adjacency matrix, defines 
the {\bf simplex graph} $(\G,\V)$ of $(G,V)$. The orientation-dependent matrix $D$ is of the form $D=d+d^*$, 
where $d$ is a lower triangular $v \times v$ matrix satisfying $d^2=0$. 
The matrix $d$ is called {\bf exterior derivative} and decomposes into blocks $d_k$ called {\bf signed incidence matrices}
in graph theory. The operator $d_0$ is the {\bf gradient}, $d_1$ is the {\bf curl}, and the adjoint $d_0^*$ is the {\bf divergence}
and ${\rm div}({\rm grad}(f)) = d_0^* d_0 f$ restricted to functions on vertices is the {\bf scalar Laplacian} $L_0=B-A$,
where $B$ is the diagonal degree matrix and $A$ is the adjacency matrix. The {\bf Laplace-Beltrami} operator 
$D^2=L$ is a direct sum of matrices $L_k$ defined on $\Omega_k$, the first block being 
the scalar Laplacian $L_0$ acting on scalar functions. To prove $d \circ d=0$ in general, take
$f \in \Omega_, x \in \G_{k+2}$ and compute $d \circ d f(x) = d( \sum_{i} f(x_i) \sigma(x_i) ) = \sum_{i,j} \sigma(i) \sigma(j(i))) f(x_{ij})$,
where $x_{i}$ is the simplex with vertex $i$ removed. Now, independent of the
chosen orientations, we have $\sigma(i) \sigma(j(i)) = -\sigma(j) \sigma(i(j))$ because permutation signatures change when the
order is switched. For example, if $x$ is a tetrahedron, then the $x_{ij}$ is an edge and each appears twice with opposite orientation.
This just given definition of $d$ could be done using antisymmetric functions $f(x_0, \dots, x_n)$, where
$df(x_0,...,x_{n+1}) = \sum_{i=0}^{n+1} (-1)^{i} f(x_0, \dots \hat{x}_i, \dots, x_{n+1})$ but the use of ordinary functions on simplices
is more convenient to implement, comes historically first and works for any finite simple graph.
A basis for all cohomology groups of any graph is obtained less than two dozen lines of Mathematica code without referring to 
any external libraries. Instead of implementing a tensor algebra, we can use linear algebra and graph data structures 
which are already hardwired into modern software. Some source code is at the end.
Dirac operators obtained from different orientations are unitary equivalent - the conjugation is given by diagonal $\pm 1$ 
matrices, even so choosing an orientation is like picking 
a basis in any linear algebra setting or selecting a specific gauge in physics. Choosing an orientation fixes a
``spin" value at each vertex of $\G$ and does not influence the spectrum of $D$,
nor the entries of $L$. Like spin or color in particle physics, it is irrelevant for
quantum $\dot{u}=iLu$ or wave mechanics $\ddot{u}=-Lu$ and only relative orientations matter for cohomology.
While the just described notions in graph theory are close to the notions in the continuum, the discrete
setup is much simpler in comparison: while in the continuum, the Dirac operator must be realized as a matrix-valued 
differential operator, it is in the discrete a signed adjacency matrix of a simplex graph.

\section{Calculus}
Functions on vertices $V=\G_0$ form the $v_0$-dimensional vector space $\Omega_0$ of {\bf scalar functions}. 
Functions on edges $D=\G_1$ define a $v_1$-dimensional vector space $\Omega_1$ of $1$-forms, functions on triangles form a 
$v_2$-dimensional vector space $\Omega_2$ etc. The matrix $d_0: \Omega_0 \to \Omega_1$ is the {\bf gradient}.
As in the continuum, we can think of $d_0f$ as a {\bf vector field}, a function on the edges given by 
$d_0( f(a,b)) = f(b)-f(a)$, think of the value of $d_0 f$ on an edge $e$ as a {\bf directional derivative} in the direction
$e$. The orientation on the edges tells us in which direction this is understood. 
If a vertex $x$ is fixed and an orientation is chosen on the edges connected to $x$ so that
all arrows point away from $x$, then $df(x,y) = f(y)-f(x)$, an identity which explains a bit more the name ``gradient".
The matrix $d_1: \Omega_1 \to \Omega_2$ maps a function on edges to a function 
on triangles. It is called the {\bf curl} because it sums up the values of $f$ along the boundary of the 
triangle and is a discrete line integral. As in the continuum, the ``velocity orientations" chosen on the edges 
and the orientation of the triangle play a role for line integrals. We check that ${\rm curl}({\rm grad}(f)) = 0$
because for every triangle the result is a sum of 6 values which pairwise cancel.
Now look at the adjoint matrices $d_0^*$ and $d_1^*$. The signed incidence matrix $d_0^*$ maps a function on edges 
to a function on vertices. Since it looks like a discretization of {\bf divergence} - we count up what ``goes in" or
``goes out" from a vertex -  we write $d_0^*f = {\rm div}(f)$. 
We have $d_0^* d_0 = {\rm div}({\rm grad}(f)) = \Delta f = L_0f$. The Laplacian $L_1$ acting on 
one-forms would correspond to $\Delta F$ for vector fields $F$, which comes close.
We have seen already that the matrix $D^2$ restricted to $\Omega_0$ is the 
{\bf combinatorial Laplacian} $L_0 = B-A$, where $B$ is the diagonal degree matrix and $A$ is the adjacency matrix. 
The story of multi-variable calculus on graphs is quickly told:
if $C$ is a curve in the graph and $F$ a vector field, a function on the oriented
edges, then the {\bf line integral} $\int_C F \; dr$ is a finite sum $\sum_e F(e)$ where $e'=1$ 
if the path goes with the orientation and $e'=-1$ otherwise. The function $e'$ on the edges of the path is the
analogue of the {\bf velocity vector}. The identity $\int_C \nabla f \; dr = f(b)-f(a)$ is due to cancellations inside
the path and called the {\bf fundamental theorem of line integrals}. A {\bf surface} $S$ is a collection of triangles in $G$ for which 
the boundary $\delta S$ is a graph. This requires that orientations on intersections of triangles cancel; otherwise 
we would get a chain. Abstract simplicial complexes is the frame work in which the story is usually told but this structure
needs some mathematical maturity to be appreciated, while restricting to calculus on graphs does not. 
Given a function $F \in \Omega^2$, define the flux $\int_S F \; dS$ as the finite sum 
$\sum_{t \in S} F(t) \sigma(t)$, where $\sigma(t)$ is the sign of the orientation of the simplex $t$. 
While $\delta S$ is a {\bf chain} in general, element in the free group generated by $\G$, the assumption that $\delta S$ is a graph
leads immediately to {\bf Stokes theorem} $\int_S F \; dS = \int_{\delta S} C \; dr$. 
If the surface $S$ has no boundary, then the sum is zero. An example is when $S$ is a graph which is a triangularization of a compact orientable 
surface and the triangles inherit the orientation from the continuum. 
The exterior derivative $d_k$ is the block part of $D$ mapping $\Omega_k$ to $\Omega_{k+1}$. It satisfies $d_{k+1} d_k = 0$.
Stokes theorem in general assures that $\int_R df = \int_{\delta R} f$ if $f$ is in $\Omega_{k-1}$ and $R$ is a subset in $\G_k$ for which 
$\delta R$ is a graph. For $k=3$, where $R$ is a collection of $K_4$ subgraphs for which the boundary $S=\delta R$ is a graph,
the later has no boundary and $\int_R df = \int_S f$ is the {\bf divergence theorem}. We have demonstrated in this paragraph that 
tools built by Poincar\'e allow to reduce calculus on graphs to linear algebra. In general, replacing continuum structures by 
finite discrete ones can reduce some analysis to combinatorics, and some topology to graph theory.

\section{Cohomology}
The $v$-dimensional vector space $\Omega$ decomposes into subspaces 
$\oplus_{k=0}^{n} \Omega_k$, where $\Omega_k$ is the vector space of $k$-forms, functions on the set
$\G_k$ of $k$-dimensional simplices in $G$. 
If the dimension of $\Omega_k$ is $v_k$, then $v = \sum_{k=0} v_k$ and $v_k$ is the number of
$k$-dimensional cliques in $G$. The matrices $d_k$ which appear as lower side diagonal blocks in $D$ belong to linear maps 
$\Omega_k \to \Omega_{k+1}$. The identity $d_k \circ d_{k-1}=0$ for all $k$ is equivalent to the fact that $L=D^2$ has a 
block diagonal structure. The kernel $C_k$ of $d_k: \Omega_k \to \Omega_{k+1}$ is called the vector 
space of {\bf cocycles}. It contains the image $Z_{k-1}$ of $d:\Omega_{k-1} \to \Omega_k$, the vector space of 
{\bf coboundaries}. The quotient space $H^k(G)=C_k(G)/Z_{k-1}(G)$ is called the {\bf $k$'th cohomology group} of the graph. Its
dimension $b_k$ is called the $k$'th {\bf Betti number}. With the {\bf clique polynomial} $v(x) = \sum_{k=0} v_k x^k$ 
and the {\bf Poincar\'e polynomial} $p(x)= \sum_{k=0} b_k x^k$, the {\bf Euler-Poincar\'e formula} tells 
that $v(-1)-p(-1)=0$, we have $v(1)=v$ and $p(1) = {\rm dim}({\rm ker}(D))$ so that $v(1)-p(1)$ 
is the number of nonzero eigenvalues of $D$. Linear algebra relates the topology of the graph 
with the matrix $D$. The proof of the Euler-Poincar\'e formula needs the {\bf rank-nullity theorem}
$v_k = z_k + r_k$, where $z_k={\rm ker}(d_k)$ and $r_k={\rm im}(d_k)$ as well as the formula
$b_k = z_k-r_{k-1}$ which follows from the definition $H^k(G)=C_k(G)/Z_{k-1}(G)$: the dimension of a quotient
vector space is equal to the dimension of the orthogonal complement. If we add up these two equations,
we get $v_k-v_k=r_k-r_{k-1}$. The sum over $k$, using $r_0=r_{n+1}=0$, telescopes to $0$. 
We have proven the {\bf Euler-Poincar\'e formula}. 
For example, if $G$ has no tetrahedra, then the Euler-Poincar\'e formula tells for a connected graph that 
$v-e+f=v_0-v_1+v_2= 1-b_1+b_2$. If $G$ is a triangularization of a surface of genus $g$, then $b_2=1,b_1=2g$ 
which leads to the formula $v-e+f=2-2g$. An other example is the circular graph $G=C_n$, where $v_0=v_1=n$ and $b_0=b_1=1$. 
A constant nonzero function on the vertices is a representative of $H^0(G)$ and a constant nonzero function
on the edges is a representative of $H^1(G)$. If $H^0(G)$ and $H^n(G)$ are both one-dimensional and all other
$H^k(G)$ are trivial, the graph is called a {\bf homology sphere graph}. It can be a triangularization of the usual sphere but
does not need to be as a triangularization of the Poincar\'e homology sphere. For 
triangularizations of compact $n$-dimensional manifolds, we have the {\bf Poincar\'e duality}
$b_k=b_{n-k}$, the reason being that the dual triangularization has the same Poincar\'e polynomial and that the dual triangularization is homotopic,
leading to the same cohomologies. This can be generalized to graph theory: whenever $G$ has the property that there exists a dual graph which is 
homotopic, then Poincar\'e duality holds. The Betti numbers $b_k$ are significant because they are independent of homotopy deformations of the graph 
as we will see below. We will also see that they can be read off as $\dim(\ker(L_k))$, so that everything is reduced to linear algebra. 
These ideas were all known to Poincar\'e, who used finite-dimensional combinatorial theory to understand what was
later called {\bf de Rham cohomology} of differential forms on compact smooth manifolds. 
Historically, the setup in the continuum needed more work and Georges de Rham was among the mathematicians who made the connection 
between the discrete and the continuum rigorous. This paragraph should have shown that the cohomology theory on graphs 
just needs some knowledge of linear algebra.

\section{Super symmetry}
The vector space $\Omega$ is the direct sum of the {\bf Bosonic subspace}
$\Omega_b = \oplus \Omega_{2k}$ and the {\bf Fermionic subspace}, the span
of all $\Omega_{2k+1}$. Let $P$ be the diagonal matrix which is equal to $1$ on $\Omega_b$ and equal to $-1$
on $\Omega_f$. For any linear map on $A$, define the {\bf super trace} $\str(A) = \tr(A P)$. 
The relations $L=D^2, 1=P^2, DP+PD=0$ define what Witten called {\bf super symmetry} in 0 dimensions.
We prove here two symmetries of the spectrum of $D$: the first is the symmetry $\lambda \leftrightarrow -\lambda$. 
The spectrum of $D$ is symmetric with respect to $0$ because if $\lambda$ is an eigenvalue of $D$
then $-\lambda$ is an eigenvalue: if $Df=\lambda f$, then $PD Pf =-DP Pf = -D f = -\lambda f$.
The vector $P f$ is the "antiparticle" to the "particle" $f$. Apply $P$ again to this identity to 
get $D (Pf) = -\lambda (Pf)$. Therefore, $Pf$ is an eigenvector to the eigenvalue $-\lambda$.
Now, we show that if $\lambda$ is an eigenvalue of $L$ to the eigenvector $f$ then it is also an eigenvalue
of $L$ to the eigenvector $Df$. If $Lf = \lambda f$, then $L (Df) = D Lf = D \lambda f = \lambda (Df)$. 
The two eigenvectors $f,Df$ belonging to the eigenvalue $\lambda $ are perpendicular if we have chosen a basis so that
$\Omega_f$ and $\Omega_b$ are perpendicular. The nonzero eigenvalues in $\Omega_b$ therefore pair with nonzero eigenvalues in 
$\Omega_f$. This leads to the identity $\str(L^k)=0$ for all $k>0$ and to 
the {\bf McKean-Singer formula} $\str(e^{-t L}) = \chi(G)$. More generally, $\str(\exp(f(D))= \chi(G)$ 
for any analytic function which satisfies $f(0)=0$. The reason is that the even part of $\exp(f(D))$ can be written as $g(L)$
and that the odd part is zero anyway because there are only zeros in the diagonal. By continuity and Weierstrass, 
it is even true for all continuous functions. 
To summarize, we have seen that the spectrum $\sigma(D)$ of $D$ is determined by the spectrum of $L$ and given as $\pm \sqrt{\sigma(L)}$
as a set. Because every nonzero eigenvalue of $D$ comes as pair $-\lambda,\lambda$, each nonzero eigenvalue of $L$ had
to appear twice. The McKean-Singer symmetry gave even more information: the pairs belong to different subspaces 
$\Omega_f$ and $\Omega_b$. Let us mention more anti-commutation structure which goes under the name ``super symmetry".
Let $\{A,B \}=AB+BA$ denote the {\bf anti-commutator} of two square matrices $A,B$. We have $\{ d_k,d_l \; \} = \delta_{lk} L_k$, where
$\delta_{lk}$ is the Kronecker delta. While this again reflects the fact that $D^2=L$ is block diagonal, it explains a bit
the relation with the {\bf gamma matrices} $\gamma^k$ introduced in the continuum to realize the anti-commutation 
relations $\{ \gamma^k, \gamma^l \; \} = -2 \delta_{lk}$ and which was the starting point of Dirac to factor the 
d'Alembert Laplacian $-\Delta$. The first order differential operator $D= \sum_i \gamma^i \delta_i$ is now the square root
of $L$. Such gymnastics is not necessary in the discrete because the exterior derivative and exterior bundle formalism
are already built into the graph $\G$. 

\section{Hodge theory} 
Vectors in the kernel of $L_k$ are called {\bf harmonic $k$-forms}. For graphs, where $L$
is a finite matrix, the Hodge theory is part of linear algebra. Our goal is to prove the {\bf Hodge theorem}, which states
that the dimension ${\rm dim}({\rm ker}(L_k))$ of the harmonic $k$-forms is equal to $b_k$. The theorem also assures that
cohomology classes are represented by harmonic forms.
Because of $\langle d_0 f,d_0 f \rangle = \langle f,Lf \rangle$, the kernel of $L_0= d_0^* d_0$ is the 
same as the kernel of $d_0$, which consists of all functions $f$ which are locally constant. The number 
$b_0$ therefore has an interpretation as the number of connected components of $G$. But lets start with the proof:
first of all, we know that $Lf=0$ is equivalent to the intersection of $df=0$ and $d^*f=0$ because $\langle f,Lf \rangle
= \langle df,df \rangle + \langle d^*f,d^*f \rangle$ shows that $Lf=0$ is equivalent to $df=d^*f=0$. 
We know already that the kernel of $L$ and the image of $L$ are perpendicular because $L$ is a symmetric matrix. 
We also know from $\langle dg,d^*h \rangle = \langle d^2 f,h \rangle = 0$ that the image of $d$ and the image of $d^*$ are 
perpendicular. They obviously together span the image of $D$ and so the image of $L=D^2$. We therefore have an orthogonal decomposition 
$\R^v = {\rm im}(d) + {\rm im}(d^*) + {\rm ker}(L)$ which is called the {\bf Hodge decomposition}. 
We use this decomposition to see that any equivalence class $[f]$ of a cocycle $f$ satisfying $df=0$ 
can be associated with a unique harmonic form $h$, proving so the {\bf Hodge theorem}. 
To do so, split $f=du+d^*v + h$ into a vector $du$ in ${\rm im}(d)$ a vector $d^*v$ in ${\rm im}(d^*)$ and 
a vector $h$ in the kernel of $L$. 
We see that $f$ is in the same cohomology class than $d^*v + h$. But since $df=0$ and $dh=0$ imply $d d^* v=0$
and because of $d^* d^*v=0$, also $d^*v$ must be in $\ker(H)$. Being in the kernel while also being perpendicular 
to the kernel forces $d^*v$ to be zero. This means that the equivalence classes are the same $[f]=[h]$. 
The just constructed linear map $\phi: [f] \to h$ is injective because a nontrivial 
kernel vector $f$ satisfying $\phi(f)=0$ would mean $f=du$ and so $[f]=0$ in $H^k(G)$. Hodge theory is useful: we not
only have obtained the dimension $b_k(G)$ of the $k$'th cohomology group, a basis for the kernel of $L_k$ gives us a concrete
basis of $H^k(M)$. Furthermore, the {\bf heat kernel} $e^{-t L}$ converges in the limit $t \to \infty$ to a 
projection $K$ onto the kernel of $L$, the reason being that on a subspace to a positive eigenvalue the heat kernel decays exponentially. 
We also see the connection with Euler-Poincar\'e: for $t=0$, the heat kernel is
the identity which has the super trace $v_0-v_1+v_2- \cdots$. In the limit $t \to \infty$, we get the super trace
of the projection $K$ which is equal to $b_0-b_1 + \dots $. 

\section{Perturbation theory}
If $G_1,G_2$ are two graphs on the same set of vertices, we can compare the spectra of their
Dirac matrices or Laplacians. Comparing the Dirac operators is easier because the entries only take the values
$1,-1$ or $0$. A nice tool to study perturbations is {\bf Lidskii theorem} which assures that if $A,B$ 
are symmetric with eigenvalues $\alpha_1 \leq \alpha_2 \leq \dots \leq \alpha_n$ and
$\beta_1 \leq \beta_2 \leq \dots \leq \beta_n$, then 
$(1/n) \sum_{j=1}^n |\alpha_j - \beta_j| \leq (1/n) \sum_{i,j=1}^{n} |A-B|_{ij}$.
The left hand side is a {\bf spectral distance}. We have rescaled so that $d(A,0) = \tr(|A|)/n$, 
makes sense for {\bf graph limits} obtained from finer and finer triangularizations of 
manifolds. Define the {\bf simplex distance} $d(G,H)$ of two graphs $G,H$ with vertex 
set $V$ as $(1/v)$ times the number of simplices of $G,H$ which are different 
in the complete graph with vertex set $V$. If simplex degree of $x \in \G_k$ is defined as the
number of matrix entries in the column $D_x$ which are not zero. We immediately get 
that $d(\sigma(G),\sigma(H)) \leq \deg \cdot d(G,H)$.
I learned the above inequality from Yoram Last who 
reduced it to a theorem of Lidskii, which deals with the sum of the eigenvalues 
$\gamma_i$ of the difference $C=A-B$ of two selfadjoint matrices $A,B$. 
Here is Yoram's argument: if $U$ is the orthogonal matrix diagonalizing $C$, then 
$\sum_j |\alpha_j-\beta_j|  \leq \sum_i |\gamma_i|$ $=\sum_i (-1)^{m_i} \gamma_i$ for some integers $m_i$.
This is equal to $=\sum_{i,k,l} (-1)^{m_i} U_{ik} C_{kl} U_{il}$
$\leq \sum_{k,l} |C_{kl}| \cdot |\sum_i (-1)^{m_i} U_{ik} U_{il}|$ $\leq \sum_{k,l} |C_{kl}|$. 
The statement has so be reduced to the classical Lidskii theorem $\sum_j |\alpha_j-\beta_j|  \leq \sum_i |\gamma_i|$
whose proof can be found in Barry Simon's trace ideal book. Having been a student in Barry's trace ideal 
course he gave in 1995, I should be able to give a proof of the later too: it uses that a matrix $A(t)$ which 
depends on a parameter $t$ satisfies the {\bf Rayleigh formula} $\lambda'=\langle u,A'u \rangle$ if $u$ is 
a normalized eigenfunction. Proof: differentiate $Au=\lambda u$ to get $A'u+Au'=\lambda' u + \lambda u'$ and 
take the dot product of this identity with $u$, using that $|u|^2=1$ gives
$\langle u,u'\rangle=0$ and that the symmetry of $A$ implies 
$\langle u,A u' \rangle=\langle Au,u'\rangle = \lambda \langle u,u' \rangle=0$. 
Using Rayleigh, we can now estimate how an eigenvalue $\lambda(t)$ of $A(t)$ 
changes along the path $A(t) = A+tC$ with constant velocity $A'(t)=C$ 
connecting $A$ with $B$. If $u(t)$ is a normalized eigenfunction of $A(t)$,
then $\lambda'(t) = \langle u(t),C u(t) \rangle$. By writing $u(t) = \sum_k a_k(t) v_k$, where $v_k$
is the eigen basis of $C$, we have $\sum_{k=1}^n a_k^2=1$. We see that $\lambda'(t) = \sum_k a_k^2 \gamma_k$
with $\sum a_k^2=1$ implying the existence of a path of stochastic matrices $S(t)$ such that 
$\lambda'_i(t) = \sum_j S_{ij}(t) \gamma_j$. Integrating this shows that $\alpha_i-\beta_i = \sum_j S_{ij} \gamma_j$
where $S$ is {\bf doubly sub stochastic} meaning that all row or column sums are $\leq 1$. Since such a matrix
is a contraction in the $l_1$ norm, $\sum_j |\alpha_j-\beta_j|  \leq \sum_i |\gamma_i|$ follows. 

\section{Cospectral graphs}
The McKean-Singer result assures that some isospectral graphs for the graph Laplacian
inherit the isospectral property for the Dirac operator. The question of Dirac isospectrality
is also studied well in the continuum. Berger mentions in his ``Panorama" that Milnor's examples 
of isospectral tori already provided examples of manifolds which are isospectral with respect to 
the Hodge Laplacian $L$. Finding and studying isospectral manifolds and graphs has become an industry. 
One reason why the story appears similar both in the continuum and discrete could be that typically
isospectral sets of compact Riemannian manifolds are discrete and for hyperbolic manifolds it is
always impossible to deform continuously. 
We can use the McKean-Singer symmetry to see why it is not uncommon that an isospectral 
graph for $L_0$ is also isospectral for $L$. This is especially convenient, if the graph has
no $K_4$ subgraphs and only a few isolated triangles. This assures that $L_2$ are isospectral. By 
McKean-Singer, then also $L_1$ has to be Dirac isospectral. An other example is given by trees 
which have no triangles so that $L=L_0 + L_1$. Therefore, if a tree is isospectral with respect to 
$L_0$, it is also Dirac isospectral. 
There are various other examples of Dirac isospectral graphs which are not graph isomorphic. 
The simplest examples are {\bf Dehn twisted flat tori}: take a triangularization of a two-dimensional torus, 
cut it along a circle, twist it along the circle, then glue them together again. Again, the reason is 
super symmetry: there is isospectrality with respect to $L_0$ and because $L_2$ is the same, also 
the $L_1$ are isospectral by McKean Singer. 
A concrete pair of graphs which is isospectral with respect to the Dirac operator 
is an example given by Hungerb\"uhler and Halbeisen. Their examples are especially interesting because it was not
brute force search, but an adaptation of a method of Gordon in the continuum which led to these graphs. 
There are also examples of graphs which are isospectral with respect to $L_0$ but not with respect to $L_1$. 
Again, this can be seen conveniently without any computation from McKean-Singer, if one has a pair for which $L_2$ are
not isospectral, then also $L_1$ is not isospectral. Most open problems from the Laplacian go over to the Dirac case: 
how much information about the geometry of $G$ can be extracted from the spectrum of the Dirac operator? 
To which degree does the spectrum determine the graph? Are Dirac isospectral
graphs necessarily homotopic in the sense of Ivashchenko discussed below? The examples we have seen
are. While the Dirac spectrum determines the Betti numbers by counting zero eigenvalues, one can ask to which degree
the traces or more generally the values of the zeta function determine the topology. 

\section{Homotopy deformations}
Lets call a graph with one vertex {\bf contractible}. Inductively, lets call a graph $G$ 
{\bf simply contractible} if there is a vertex $x$ of $G$ such that its unit sphere $S(x)$ is contractible and $G$ without 
$x$ and connecting edges is simply contractible.  A {\bf homotopy expansion} is an addition of a 
new vertex $x$ and connections to $G$ such that in the sub graph $S(x)$ in the new graph 
is contractible. A homotopy reduction is the removal 
of a vertex and all its edges for which $S(x)$ is contractible. Two graphs are called {\bf homotopic} if there is 
a sequence of homotopy steps which brings one into the other. A graph is called {\bf contractible} if it is homotopic
to a one vertex graph.  These notion shadow the definitions
in the continuum and indeed, Ivashchenko who defined this first, was motivated by analogous notions put
forward by Whitehead. The definitions since been simplified by Chen-Yau-Yeh. 
Discrete homotopy is simple and concrete and can be implemented on a computer. It is useful theoretically and also in applications
to reduce the complexity to calculate things. 
It motivates the notion of ``critical points" and ``category" and ``index" of a critical point. {\bf Morse theory}
and {\bf Lusternik-Schnirelmann theory} go over almost verbatim. For a function $f$ on the vertex set $V$, a vertex 
$x$ is called a {\bf critical point} if $S^-(x)=S(x) \cap \{ y \; | \;  f(y)<f(x) \; \}$ is not contractible. 
The index of a critical point is defined as $i_f(x)=1-\chi(S^-(x))$, where $\chi$ is the Euler characteristic. As in the
continuum there can be critical points which have zero index. 
By induction, building a graph up one by one one can see that the {\bf Poincar\'e-Hopf theorem} $\sum_x i_f(x)=\chi(G)$ 
holds. This is completely analogue to Morse theory, where manifolds change by adding handles when a critical point is
reached by a gradient flow of a Morse function $f$. 
When averaging the index $i_f(x)$ over functions we obtain curvature $K(x)$. Curvature is defined for any simple graph
and defined as follows: if $V_k(x)$ the number of subgraphs $K_{k+1}$ of the unit sphere $S(x)$ with the assumption
$V_{-1} = 1$ then $K(x)$ is defined as
$\sum_{k=0}^{\infty} (-1)^k V_{k-1}(x)$. The Gauss-Bonnet formula $\sum_x K(x) = \chi(G)$ is a direct consequence of the 
Euler handshaking lemma $\sum_{x \in V} V_{k-1}(x) = (k+1) v_k$. 
Curvature becomes more geometric if graphs have a geometric structure. While Gauss-Bonnet is surprisingly simple, it 
is also astonishingly hard to prove that $K$ is identically zero for odd dimensional graphs, where every sphere $S(x)$ has
Euler characteristic $2$. The vector space $\Omega$ has also an algebra structure. This exterior
product produces the {\bf cup product} on cohomology. There is a relation between {\bf cup length} which is 
an algebraic notion, and {\bf Lusternik-Schnirelmann category} and the minimal number of critical points a function $f$ on
a graph homotopic to $G$ can have. The relation between Morse theory, homotopy, category and 
the algebra appear here in a finite setting. Homotopy is an equivalence
relation on graphs. It is probably an impossible task to find a formula for 
the number  homotopy types there are on graphs with $n$ vertices. While most of these notions are known 
in one way or the other in the continuum, they are much more complicated in the later. The index averaging
result ${\rm E}[i_f(x)] = K(x)$ has no natural continuum analogue yet, notably because we have no natural probability space of 
Morse functions on a manifold yet. For Riemannian manifolds $M$ embedded in an ambient Euclidean space (which is always possible
by Nash), one can use a finite dimensional space of linear functions in the ambient space to induce functions on $M$ 
and see the differential geometric curvature as the expectation of the index $i_f(x)$ of functions. This provides a
natural integral geometric proof of the Gauss-Bonnet-Chern for compact Riemannian manifolds. Seeing curvature as expectation
has other advantages like that if we deform a structure, then we can push forward the measure on functions and get new 
curvature without having to drag along the Riemannian formalism. Different measures define different curvatures for which
Gauss-Bonnet holds automatically. 

\section{Simplex combinatorics}
Let $\deg_p(x)$ be the number $(p+1)$-dimensional simplices which contain the $p$-dimensional simplex $x$. 
This generalizes $\deg_0(x)$, which is the usual degree of a vertex $x$. For $p>0$ we have ${\rm deg}_p(x) = L_p(x,x)-(p+1)$,
the reason being that $L_p(x,x)$ counts also the $p+1$ connections with the $p+1$, $p-1$ dimensional simplices in 
the simplex $x$. For example, if $x$ is a triangle, where $p=2$, then there are three one-dimensional sub-simplices given by edges 
inside $x$ which are connected to $x$ in $\G$. We also have ${\rm deg}_0(x) = L_0(x,x)$. It follows that $\tr(L_p) = (p+2) v_{p+1}$. 
A {\bf path} in $\G = \bigcup \G_k$ is a sequence of simplices
contained in $\G_k \cup \G_{k+1}$ or $\G_k \cup \G_{k-1}$ such that no successive $i=x_k,j=x_{k+1}$ have the property that $D_{ij} \neq 0$.
In other words, it is a path in the graph $\G$. If $k$ is even, then the path is called even, otherwise odd. 
The integer $|D|^k_{xy}$ is the number of paths of length $k$ in $\G$ starting at a simplex $x$
and ending at a simplex $y$. The entry $L^{k}_{xx} = [D^{2k}]_{xx}$ is the number of closed paths of length $2n$. 
Note that $D^{2k+1}_{xx}=0$ which reflects the fact that any step changes dimension and that the start
and end dimension is the same. The trace $\tr(L^k) = \tr(D^{2k})$ is therefore total number of closed paths in
$\G$ which have length $2k$. Because $\str(L^k)=0$, we get as a
combinatorial consequence of the McKean-Singer symmetry that the number of odd paths is the same than the 
number of even paths. Since for any graph $G$, the simplex graph $\G$ is a new graph, one could ask what the simplex graph 
of $\G$ looks like. But this is not interesting since $\G$ has no triangles. For a triangle $G=K_3$ 
for example, $\G$ is the cube graph with one vertex removed. Its simplex graph is a subdivision of it
where every edge has got a new vertex. Repeating that leads to more and more subdivisions and all these
graphs are what one calls homomorphic. They are also homotopic since in general, homomorphic graphs are
homotopic. By the Kirchhoff matrix tree theorem applied to $\G$, the pseudo determinant of $B-|D|$ divided
by $v$ is the number of spanning trees in $\G$, where $B$ is the diagonal degree matrix in $\G$. 
There is a better matrix tree theorem in $\G$ which relates the pseudo determinant of $L=D^2$ with the number
of trees in $\G$. We will look at it below. 

\section{Variational questions}
Since graphs are finite metric spaces, one can look at various quantities and try to 
extremize them. One can look for example at the {\bf complexity} $\Det(L_0)$ which is $n$ times the
number of spanning trees in $G$, or the Euler characteristic $\chi(G) = \str( \exp(-L))$ for which 
the question of minimal $\chi(G)$ looks interesting.  Well studied is the
{\bf smallest possible rank}, a matrix with nonzero entries at nonzero entries of the adjacency matrix
can take. Recently, one has looked at the {\bf magnitude} $|G|= \sum_{i,j} Z^{-1}_{ij}$, defined by 
the matrix $Z_{ij} = \exp(-d(i,j))$, where $d(i,j)$ is the geodesic distance between to vertices $i$ and $j$. 
This quantity was forward by Solow and Polasky. We numerically see that 
among all connected graphs, the complete graph has minimal magnitude and the star graph maximal magnitude. 
The magnitude defined for any metric space so that one can look for a general metric space at the 
supremum of all $|G|$ where $G$ is a finite subset with induced metric. The {\bf convex magnitude conjecture} 
of Leinster-Willington claims that for convex subsets of the plane, 
$|A|=\chi(A) + p(A)/4 - a(A)/(2\pi)$, where $p$ is the perimeter and $a$ the area. 
The sum $\sum_{i,j} \exp(-t L_0)_{ij}$ is always equal to $v_0$ but
$\sum_{i,j} \exp(-t D)_{ij}$ resembles the magnitude. One can also look at $\sum_{i,j,p} (-1)^p \exp(-t L_p)_{ij}$. 
Other variational problems are to find the graph with maximal Laplacian permanent $\per(L)$, 
the adjacency matrix permanent $\per(A(G))$ or maximal Dirac permanent $\per(D)$. The permanent has been studied and is
combinatorially interesting because for a complete graph $K_n$, the integer $\per(A(K_n))$ is the set of permutations of 
$n$ elements without fixed point. In all cases for which we have computed $\per(D)$, the later only took a few values and 
most often $0$. Also $\per(D)=4$ occurred relatively often for some non-contractible graphs. 
Since it is difficult however to compute the permanent $\per(D)$ even for smaller graphs, we have no clue
yet how a distribution of $\per(D(G))$ would look like nor what the relation is to topology.
An other quantity motivated from the continuum are notions of {\bf torsion}. 
McKean-Singer implies that the $\eta$-function $\eta(s) = \sum_{\lambda \neq 0} (-1)^k \lambda^{-s}$
is always $0$ for graphs. The number $\exp(-\eta'(0))$, the analogue of
the {\bf analytic torsion}, is therefore always equal to $1$. Analytic torsion for $D$ can be rewritten as
$\tau(G) = \exp(\str(\log(\Det(L)))) = (\prod_{\lambda \in \sigma_b} \lambda)/(\prod_{\lambda \in \sigma_f} \lambda) = 1$, 
where $\sigma_f,\sigma_b$ are the Fermionic and Bosonic eigenvalues. In the continuum, one has looked at functionals like
$(\prod_{\lambda \in \sigma_b} \lambda^{p(\lambda)})/(\prod_{\lambda \in \sigma_f} \lambda^{p(\lambda)})$, where $p(\lambda)=p$
if $\lambda$ is an eigenvalue of $L_p$. Maybe this leads to interesting functionals for graphs; interesting meaning that
it is invariant under homotopy deformations of the graph. 

\section{Matrix tree theorem}
The {\bf pseudo determinant} of a matrix $D$ is defined as 
the product of the nonzero eigenvalues of $D$. It is a measure for the complexity of the graph. 
The classical matrix tree theorem relates the pseudo determinant of the graph Laplacian $L_0$ 
with the number of spanning trees in $G$. There is an analogue combinatorial description 
for the pseudo determinant $\Det(L)$  of the Laplace-Beltrami operator $L$. The result is based on a generalization of
the {\bf Cauchy-Binet formula} for pseudo determinants $\Det(A)$. While the identity $\det(A B) = \det(A) \det(B)$
is false for pseudo determinants $\Det$, we have found that $\Det^2(A) = \sum_P \det^2(A_P)$, where the sum
is over all matrices $A_P$ of $A$ for which $b = \sum_i b_i$ redundant columns and rows have been deleted. 
The classical Kirchhoff matrix tree theorem gives an interpretation of $\Det(L_0)/n$ as the number
of spanning trees of the graph. An interpretation of $\Det(L)$ 
is given as a weighted number of spanning trees of a double cover of the simplex graph $\G$
where trees are counted negative if they have different type on the two branches and 
positive if they have the same type. This is just a combinatorial description of what $\det^2(A_P)$ means. 
Our new Cauchy-Binet result is actually more general and gives the coefficients of the characteristic polynomial 
$p(x) = \det(A-x)$ of $A=F^T G$ for two $n \times m$ matrices $F,G$ as a product of minors:
$p_k = (-1)^k \sum_{P} \det(F_P) \det(G_P)$ where $P$ runs over all possible $k \times k$ minors $P=P(IJ)$ of $F,G$. 
This is a {\bf quadratic analogue} of the well known trace formula $p_k(A) = (-1)^k \sum_P \det(A_P)$ 
of the characteristic polynomial $\det(A-x)=p_k (-x)^k + ...$ of a square matrix $A$ as the sum of 
symmetric minors $P=P(II)$ of size $k={\rm rank}(A)$. The Cauchy-Binet theorem implies that
$p_k(A^2) = (-1)^k \sum_P \det(A_P)^2$, where $P$ runs over all possible $k \times k$ minor masks $P=P(IJ)$. It implies
the {\bf Pythagoras theorem for pseudo determinants}: 
if $A=A^*$ has rank $k$, then $\Det^2(A) = \sum_P \det^2(A_P)$, where $\det(A_P)$ runs over all 
$k \times k$ minors $P=P(IJ)$ of $A$ defined by choosing row and column subsets $I,J$ of cardinality $k$. 
Having looked hard to find this Pythagoras result in the literature, we could not find it anywhere. 
Note that it is a {\bf quadratic relation} unlike the well known $\Det(A) = \sum_P \det(A_P)$, where
$A_P$ runs over all symmetric minors $P$ of ${\rm rank}(A)$ and where no matrix multiplication is involved.
Cauchy-Binet by definition deals with the matrix  product.
The proof of this general Cauchy-Binet result is quite elegant with multi-linear algebra because the trace formula
for $A=F^T G$ shows that $p_k = (-1)^k {\rm tr}(\Lambda^k F^T G) = \sum_I \sum_K \det(F_{IK}) \det(G_{KI})$
 is just a reinterpretation of what the matrix product in $\Lambda^k M(m,n) \times \Lambda^k M(n,m) \to \Lambda^k M(m,m)$ 
means. The entries of $(\Lambda^k A)_{IJ}$ of $\Lambda^k A$ are labeled by subsets $I,J \subset \{ 1, \dots, n \; \}$ 
and are given by $(\Lambda^k A) _{IJ} = \det(A_{IJ})$, where the later is the minor obtained by taking the determinant of the
intersection of the $J$ columns with $I$ rows of $A$. 
While $\per(AB) \neq \per(A) \per(B)$ in general, I learned at the Providence meeting that there
are Cauchy-Binet versions for permanents. The analogue of the pseudo determinant would be 
$(-1)^k a_k$ of the {\bf permanental characteristic polynomial} $p(x) = per(A-x) = a_k x^k + a_{k+1} x^{k+1} + ...$.
Minc's book suggests that the proof done for the exterior algebra could have an analogue when using 
a completely symmetric tensor algebra. 

\section{Differential equations}
For $A \in \Omega_1$, the 2-form $F=dA$ satisfies the {\bf Maxwell equations} $dF=0,d^*F=j$. 
If $A$ has a {\bf Coulomb gauge} meaning that $d^* A=0$, then the equation is $LA=j$ so that if $j$ is in the 
ortho-complement of the kernel of $L$ then $A$ can be found satisfying this system of linear equations. 
The equation $L A = j$ is the discrete analogue of the {\bf Poisson equation}. In the vacuum
case $j=0$, it reduces to the {\bf Laplace equation} $Lu=0$. 
The {\bf transport equation} $u_t = i D u$ is a Schr\"odinger equation. As in one dimensions, the wave 
equation $u_{tt} = -Lu$ is a superposition of two transport equations $u_t = i D u$ and $u_t = - i D u$. 
The {\bf heat equation} $u_t = -L u$ has the solution $e^{-Lt}u(0)$. The Fourier idea is to solve
it using a basis of eigenfunctions of $L$. The matrix $e^{-L t}$ is the heat kernel. It is significant
because for $t \to \infty$ it converges to a projection operator $P$ from $\Omega$ onto the linear
space of harmonic functions. The {\bf wave equation} $u_{tt} = -L u$ has the solutions 
$u(t) = \cos(Dt)u(0) + i \sin(Dt) D^{-1} u'(0)$, which works if the velocity is in the 
ortho-complement of the kernel of $L$. The name "wave mechanics" for quantum mechanics is justified
because if we write $\psi = u(0) + i D^{-1} u'(0)$, then $e^{i Dt} \psi = \psi(t)$ so that the 
wave equation for real waves is equivalent to the Schr\"odinger equation for the Dirac operator. 
All equations mentioned far are linear. A nonlinear integrable system $\dot{D}=[B,D]$
will be discussed below. Other nonlinear system are {\bf Sine-Gordon} type equations like
$Lu = \sin(u)$ for $v$ variables and for which one can wonder about the structure
of the nonzero solutions $u$. Or one could look at systems $u_{tt}+L u = \sin(u)$ for which 
it is natural to ask whether it is an integrable Hamiltonian system. 
Back to the wave equation, we can ask what the significance is of the evolution of waves on $\Omega_k$.
In classical physics, one first looks at the evolution on scalar functions. For a two-dimensional membrane,
McKean-Singer actually shows that the evolution on $\Omega_0$ determines the rest, the reason being that $L_2,L_0$
behave in the same way and that $L_1$ is the only Fermionic component which is then determined. The block $L_1$ of $L$ is 
however important in electromagnetism: the Maxwell equations $dF=0,d^*F=j$ determine the electromagnetic field $F=dA$
from the current $j$ are given by $L_1 A = j$. In the continuum, where the equations live in {\bf space-time}
and $L_1$ is a {\bf d'Alembert operator}, the equation $L_1 f =0$ already is the wave equation and a harmonic $1$-form
solving it provides a solution $F=df$ of the Maxwell equations. In the graph geometrical setup, this solution $F$ is 
a function on triangles. In a four dimensional graph, one has then on any of the $K_5$ simplices $x$, there are 
$6$ triangles attached to a fixed vertex and the function values on these triangles provide 3 electric $E_1(x),E_2(x),E_3(x)$
and 3 magnetic components $B_1(x),B_2(x),B_3(x)$. Now we can produce ``light" in a purely graph theoretical setup. 
Give a current $j$ in the orthocomplement of the kernel of $L_1$ - this can be seen as a conservation law requirement - 
then $L_1 A=j$ can be solved for $A$ providing us with a field $F=dA$. If the graph is simply connected, in which case $L_1$
is invertible, we can always solve this. The simplest four dimensional ``universe" on which we can ``switch on light" is $K_5$. 
Given any function $j$ on its edges, we can find a function $F$ on the $10$ triangles. The electromagnetic field at a
vertex $x$ is then given by the values of function on the $6$ triangles attached to $x$. But it is a small,small world! 

\section{Graph automorphisms}
A graph {\bf endomorphism} $T$ on a graph $G$ induces a linear map $T_k$ on the cohomology group 
$H^k(G)$. If $T$ is invertible it is called a {\bf graph automorphism}. The set of all automorphisms forms
a group and a classical {\bf theorem of Frucht} assures that any finite group can occur as an automorphism group of a graph.
We have looked at the question of fixed points of automorphisms in order to see whether results like the famous Brouwer fixed
point theorem can be translated to graph theory. It initially looks discouraging: simple examples like a rotation on 
a triangular graph show that there are no fixed vertices even for nice automorphisms on nice graphs which are contractible.
But if one looks at simplices as basic points, then the situation changes and all the results go over to the discrete. 
The {\bf Lefschetz number} $L(T) = \sum_k (-1)^k \tr(T_k)$ can be written as a sum
of indices $i_T(x)$ of fixed points $x$, where $(-1)^{\dim(x)} \sign(T|x)$ is the discrete analogue
of the Brouwer index. As for the Gauss-Bonnet or Poincar\'e-Hopf theorems for graphs mentioned above, 
the main challenge had been to find the right notion of index. This required a lot of searching with the help 
of the computer to check examples. 
The Lefshetz fixed point theorem $L(T) = \sum_{x} i_T(x)$ implies that if $G$ is contractible in the sense that 
$H^k(G)=0$ except for $k=0$ forcing $L(T)=1$, every automorphism has a fixed simplex. If $T$ is the identity, 
then every simplex is a fixed point and since $\tr(T_k)$ is $b_k$ in that case, the result reduces then to 
the Euler-Poincar\'e formula.  For contractible graphs, it is a discrete
analogue of the Brouwer fixed point theorem. For the proof it is crucial to look at the orbits of simplices.
Also this result shows how fundamental the concept of a simplex is in a graph. For many practical 
situations, it does not really matter whether we have a fixed vertex or a fixed simplex.
Lets look at a reflection $T$ on a circular graph $C_5$. The map $T_0$ is the identity because every
equivalence class of cocycles in $H^0(G)$ is constant. The map $T_1$ maps an equivalence class $f$ in $H^1(G)$
to $-f$. This implies that $\tr(T_1)=-1$ and $L(T) = 2$. The Lefshetz fixed point theorem assures that there
are two fixed points. Indeed, a reflection on $C_5$ always has a fixed vertex and a fixed edge. 
The fixed vertex has dimension $0$ and index $1$. The fixed edge has dimension $1$ and index $(-1)^1 \cdot (-1) =1$
because $T$ restricted to $x$ induces an odd permutation with negative signature. 
The sum of the two indices is equal to the Lefshetz number $L(T)=2$.
The zeta function $\zeta_T(z) = \exp(\sum_{n=1}^{\infty} L(T^n) \frac{z^n}{n})$ 
is a rational function with a computable product formula. It is an adaptation of the Ruelle Zeta 
function to graph theory. The Ruelle Zeta function on the other hand builds on the Artin-Mazur zeta function. 
If $\A$ is the automorphism group of $G$, we have a zeta function $\zeta(G) = \prod_{T \in \A} \zeta(T)$. 

\section{Dimension}
Since a simplex $K_{k+1}$ is $k$-dimensional, a graph contains in general 
different dimensional components. How to define dimension for a general graph?
Traditionally, algebraic topology and graph theory treat graphs as one-dimensional objects, 
ignoring the higher dimensional structure given by its sub-simplices. Indeed, standard 
topological notions of dimension either lead to the dimension $0$ or $1$ for graphs. The notion
simplicial complexes are then introduced to build higher dimensional structures. One can do
just fine in a purely graph theoretical setting. The {\bf Menger-Uryson dimension} can be 
adapted to graph and becomes meaningful.  The notion is inductive: 
by defining  $\dim(G) = \frac{1}{|G|} \sum_{x \in V} 1+\dim(S(x))$ with the assumption that $\dim(\{\})=-1$. 
A graph is called {\bf geometric} of dimension $k$, if each sphere $S(x)$ is a geometric 
graph of dimension $k-1$ and has Euler characteristic of dimension $1+(-1)^k$.
A one-dimensional geometric graph is a circular graph
$C_n$ with $n \geq 4$. While $C_3$ is two-dimensional, it is not geometric because the unit
spheres are one-dimensional of Euler characteristic $1$ and not $0$ as required. Examples
of two-dimensional geometric graphs are the icosahedron and octahedron or any 
triangularization of a smooth compact two-dimensional surface. A two-dimensional torus graph
is an example of a non-planar graph. There are planar graphs like $K_4$ which are three dimensional. 
Two dimensional geometric graphs are easier to work with than planar graphs. The four color theorem
for example is easy to prove for such graphs. Lake constructions done by taking a vertex and removing
all its edges or pyramid constructions on triangles by adding a tetrahedron above a triangle does
not change the colorability. Two dimensional graphs modified like this are rich enough to model any
country map on a two dimensional manifold. Dimension is defined for any graph and can be a fraction. 
The truncated cube for example has dimension $2/3$ because each of its vertices has dimension $2/3$. 
Dimension works naturally in a probabilistic setup because we can pretty explicitly compute the
average dimension (as we can compute the average Euler characteristic) in Erd\'os-Renyi probability
spaces. We mentioned dimension also because it is unexplored yet how much dimension information 
we can read off from the spectrum of $D$. It would probably be too much to ask to read of the dimension
from the spectrum, but there should certainly be inequalities as in the continuum. 

\section{Dirac Zeta function} 
The {\bf zeta function} of a matrix with nonnegative eigenvalues is defined as
$\zeta(s) = \sum_{\lambda \neq 0} \lambda^{-s} = e^{-s \log(\lambda)}$,
where $\lambda$ runs through all nonzero eigenvalues. As a finite sum of analytic functions, 
it is always {\bf analytic}. If $A$ has only nonnegative eigenvalues, then there is no ambiguity with $\log$
and $\zeta(s)$ is unambiguously defined. In the Dirac case, the eigenvalues take both positive
and negative values and $\lambda^{-s}$ has to be fixed. One can do that by defining $\lambda^{-s} = 
(1+e^{-i \pi s}) |\lambda|^{-s}$ if $\lambda<0$.  The Dirac zeta function of a graph is then the
zeta function of $D^2$ multiplied by $(1+e^{-2\pi i s})$. This removes some poles also in the continuum. 
The Dirac zeta function of the circle is $(1+e^{i \pi s}) \zeta_R(s)$, where $\zeta_R$ is the 
Riemann zeta function because $D= i \partial_x$ has eigenvalues $n$ with eigenfunctions $e^{-i n x}$. 
For the circular graph $C_n$, which has the eigenvalues $2 \sin(\pi k/n), k=0, \dots, n-1$, the Dirac zeta function 
of $C_n$ is $\zeta(s) = (1+e^{-i \pi s}) \sum_{k=1}^{n-1} \sin^s(\pi k/n)$. 
The {\bf pseudo determinant} of $D$ is a regularized determinant because $\Det(A) = \exp(-\zeta'(0))$. 
The later formula follows from $-d/ds \lambda^{-s} = -d/ds \exp(-s \log(\lambda)) = \log(\lambda) \exp(-s \log(\lambda))$ 
which is $\log(\lambda)$ for $s=0$. 
By definition of the zeta function, we have $\zeta(-n) = \tr(D^n)$ for positive $n$. The chosen branch is 
compatible with that and gives a zero trace if $n$ is odd. This relation and the fact that the knowledge of 
the traces determines the eigenvalues shows that knowing the zeta function is equivalent to the knowledge of the eigenvalues.
In topology, one looks at the $\eta$ function $\sum_p (-1)^p \zeta_k(p)$, where $\zeta_p$ is the zeta function of $L_p$.
The analytic torsion $\exp(-\eta'(0))$ is always equal to $1$ and only becomes more
interesting for more general operators on Dirac complexes allowing 
to distinguish manifolds which are cohomologically and homotopically identical
like a cylinder and a M\"obius strip. There should be Dirac zeta functions which enable 
to distinguish such differences. 

\section{Isospectral deformation} 
Given $D=d+d^*$ look at $B=d-d^*$. The {\bf Lax pair} $\dot{D} = [B,D]$ 
defines an isospectral deformation of the Dirac operator. As for any Lax pair,
the unitary operator $U'=BU$ has the property that it conjugates $D(t)$ with $D(0)$: 
the proof is $\frac{d}{dt} U^*(t) D(t) U(t) = U^* B^* D U + U^* D' U + U^* D B U = 0$
and the fact that for $U(0)=1$.  The deformed operator $D(t)$
is of the form $D=d+d^* + b$. We can rewrite the differential equations as
$d' = d b-b d, b' = d d^* - d^* d$.They  preserve $d^2=(d^*)^2=0$ and $L=\{d,d^* \; \}$.
We also have $\{d,b \; \}=\{d^*,b \; \}=0$ for all times.
If $df=0$ is a cocycle, and $f'=b(t)f(t)$ then $d(t) f(t)=0$ so that $f(t)$ remains a
cocycle. If $f=dg$ is a coboundary and $g'=bg$ then $f(t) = d(t) g(t)$ so that $f(t)$ remains a coboundary.
The system therefore deforms cohomology: graph or de Rham cohomology does not change if we use $d(t)$ instead
of $d$. The matrix $D(t)$ has gained some block diagonal part $b(t)$
at the places where the Laplace-Beltrami operator $L$ has its blocks. This has not prevented the
relation $D^2=L$ to remain true. The later operator $L$ does not move under the deformation because $L'=[B,L]$ 
which is zero because $B^2=-L$. What is exciting however that nevertheless the deformation changes
the geometry. The deformation defines a new exterior
derivative $d(t)$ and a new Dirac operator $C=d+d^*$. We can show that ${\rm tr}(C^2)$ is 
monotonically decreasing so that $d(t)$ converges to zero and $D(t)$ converges to a 
block diagonal matrix $b$ which has the property that $V=d^2$ agrees with the Laplacian.
The above system has a scattering nature. It can be modified by defining $B=d-d^* + i b$. 
Now, the nonlinear unitary flow $U(t)$  will asymptotically converge to the linear Dirac 
flow $e^{i b(t)}$ which also leads to a solution of the wave equation on the graph. 
We can show that the trace of $M(t) = (d(t) + d(t)^*)^2$ goes to zero 
monotonically so that $M(t)$ goes to zero monotonically. This leads to an expansion of 
space with an inflationary fast start. To show this, one proves that $\tr(M(t))$ goes to 
zero monotonically. This requires to look at the change of the eigenvalues.
There is an infinite-dimensional family of isospectral Hamiltonian deformations like that but all 
have the same features. In a Riemannian manifold setting, the analysis is infinite dimensional but
essentially the same. An interesting question is what the effect on the geometry is. 
The Dirac operator determines the metric in the graph and if we take $C(t)$ as the new
Dirac operator, then space expands. Why does the evolution take place at all. The answer
is "why not?". If a system has a symmetry then it would be a particular strange situation
if the system would not move. A rigid body in space rotates. Unless locked to an
other rock leading to friction the probability to have zero angular momentum is zero. The
isospectral deformation considered here is a {\bf Noether symmetry of quantum mechanics} 
which is to a great deal invisible because $L$ and so any wave, heat or 
Schr\"odinger dynamics is not affected. 
It only affects the Dirac operator $D$. Besides the expansion with initial inflation, the evolution
has an other interesting effect: super symmetry is present, but not visible. With the full
Dirac operator $D(t)$ which has developed a diagonal part, the super pairs $f,D(t) f$ are 
no more perpendicular for nonzero $t$. Actually, their angle is exponentially small in $t$. If $f$ was a 
Boson, then $Df$ is no more a Fermion. We do not see the super symmetry any more.
It would therefore surprise if we ever could detect super-symmetry in experiments except early
in the evolution. As mentioned before, the unitary $U(t)$ defined by the deformation pushes forward
measures on functions and so allows to define curvature in a deformed setting as the deformed 
expectation. We have only started to look yet at the question what happens with the geometry 
under deformation. 

\section{An example}
The figure shows at a graph with $7$ vertices and $9$ edges homotopic to $C_4$. Since $G$ has $2$
triangles, the Euler characteristic is $0$.
The second picture shows it with the orientation used to generate the matrix $D$.
The two triangles do not form a surface because their orientations do not match on the intersection.
Changing the orientation on the triangle $(1,2,3)$ would change the signs in the second last row
and second last column of $D$. The pseudo determinant of $D$ is $1624$.
We have $v_0=7,v_1=9, v_2=2$ and $b_0=b_1=1$ showing
that $\chi(G)=0$. The curvatures are $(K(1),...,K(6))=(2,1,-2,-4,0,0,3)/6$. By Gauss-Bonnet
they add up to $1$. With the Morse function $f(x)=x$, the point $1$ has index $1$ and
the point $6$ index $-1$ all other indices are $0$. When building up the graph starting at
$1$, all additions are homotopy steps except when adding vertex $6$ because $S(6)=\{4,5\}$
has $\chi(S(6))=2$ leading to the index $-1$. Here is the matrix $D$ with division lines drawn
between different dimension blocks for clarity:  

\begin{center}
\hspace{-5mm}
\parbox{16cm}{
\parbox{7.6cm}{\scalebox{0.28}{\includegraphics{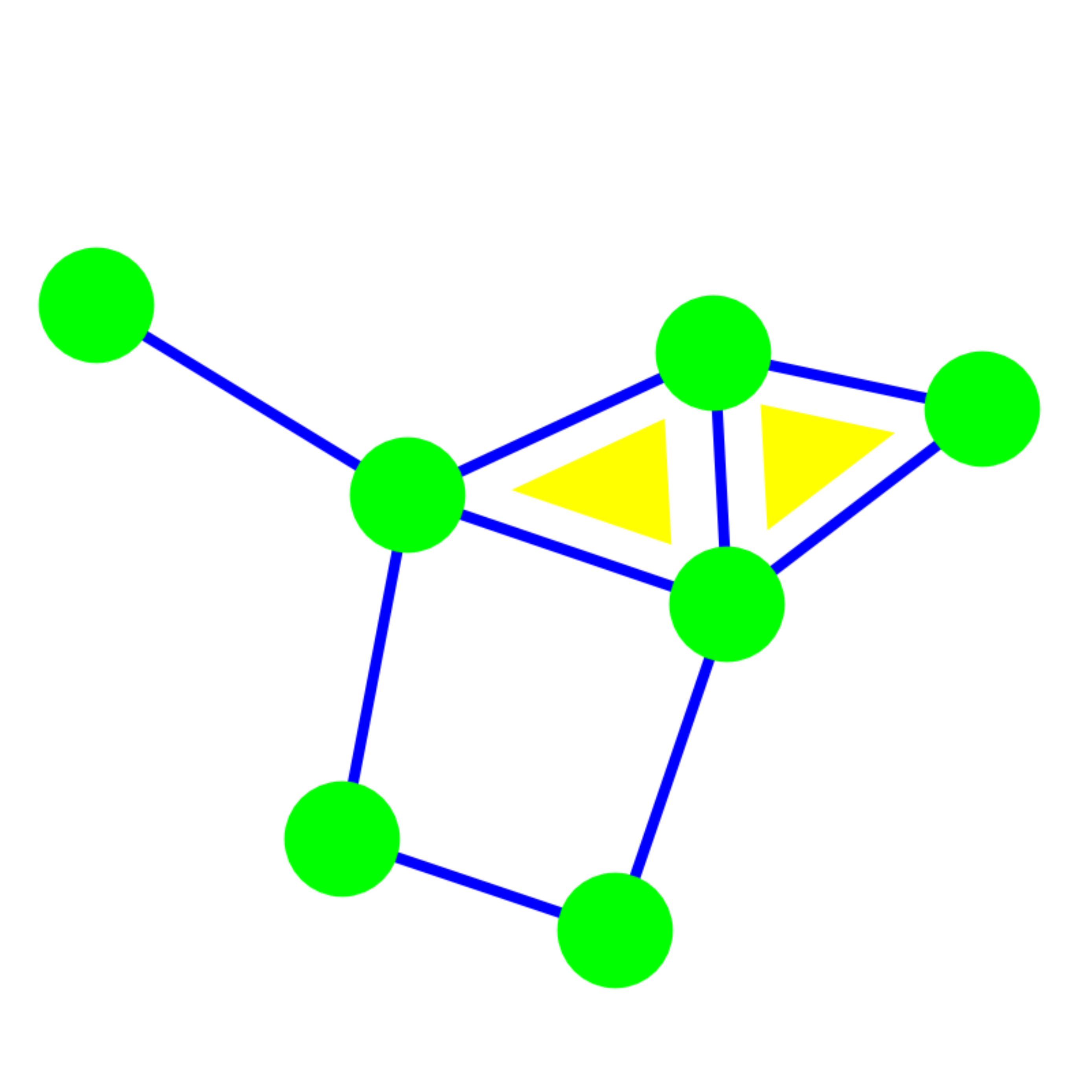}}}
\parbox{7.6cm}{\scalebox{0.28}{\includegraphics{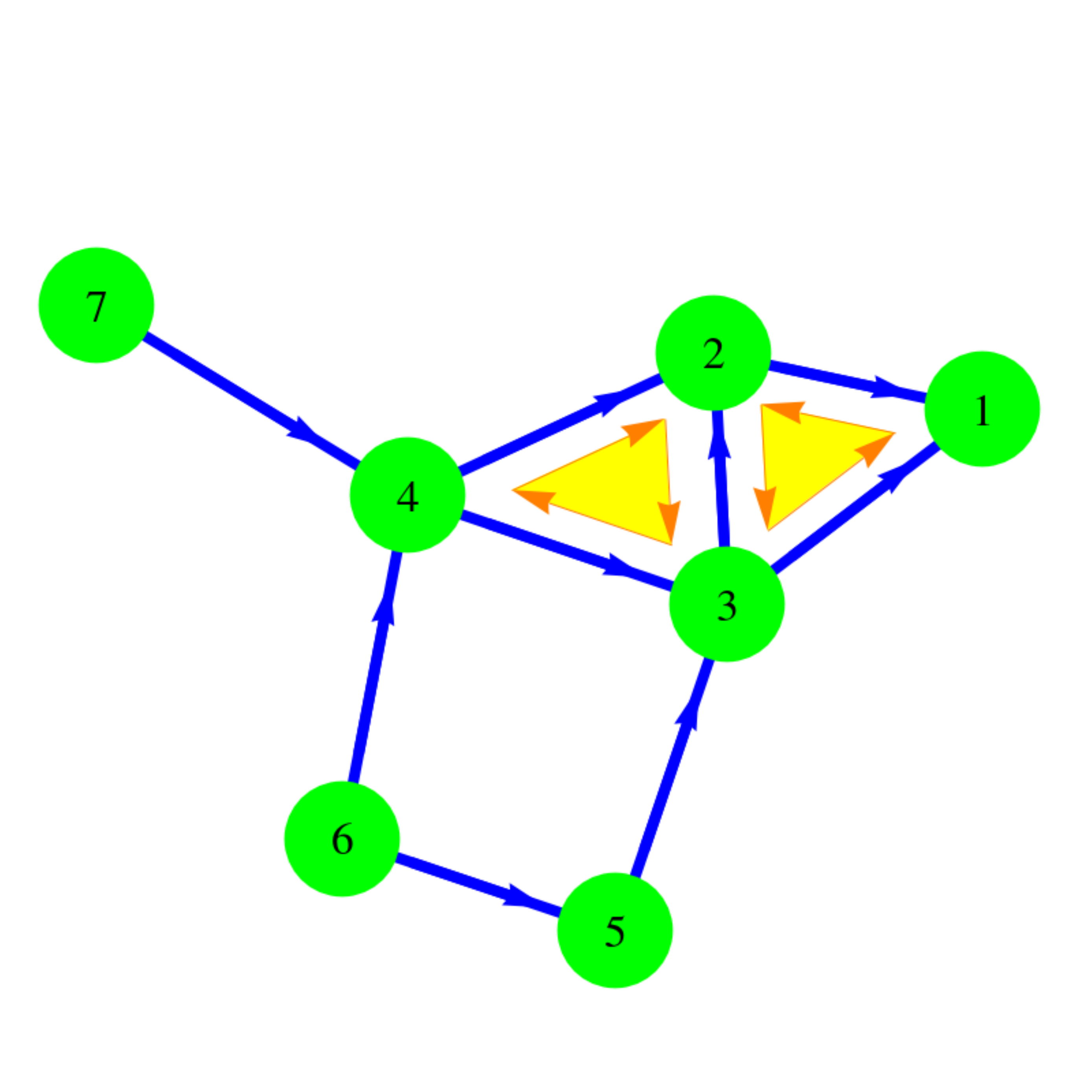}}}
}
\end{center}

\begin{tiny}
$$ D=\left[\begin{array}{ccccccc|ccccccccc|cc}
                   0 & 0 & 0 & 0 & 0 & 0 & 0 & -1 & -1 & 0 & 0 & 0 & 0 & 0 & 0 & 0 & 0 & 0 \\
                   0 & 0 & 0 & 0 & 0 & 0 & 0 & 0 & 1 & -1 & -1 & 0 & 0 & 0 & 0 & 0 & 0 & 0 \\
                   0 & 0 & 0 & 0 & 0 & 0 & 0 & 1 & 0 & 0 & 1 & -1 & -1 & 0 & 0 & 0 & 0 & 0 \\
                   0 & 0 & 0 & 0 & 0 & 0 & 0 & 0 & 0 & 1 & 0 & 1 & 0 & -1 & -1 & 0 & 0 & 0 \\
                   0 & 0 & 0 & 0 & 0 & 0 & 0 & 0 & 0 & 0 & 0 & 0 & 1 & 0 & 0 & -1 & 0 & 0 \\
                   0 & 0 & 0 & 0 & 0 & 0 & 0 & 0 & 0 & 0 & 0 & 0 & 0 & 1 & 0 & 1 & 0 & 0 \\
                   0 & 0 & 0 & 0 & 0 & 0 & 0 & 0 & 0 & 0 & 0 & 0 & 0 & 0 & 1 & 0 & 0 & 0 \\ \hline
                   -1 & 0 & 1 & 0 & 0 & 0 & 0 & 0 & 0 & 0 & 0 & 0 & 0 & 0 & 0 & 0 & 1 & 0 \\
                   -1 & 1 & 0 & 0 & 0 & 0 & 0 & 0 & 0 & 0 & 0 & 0 & 0 & 0 & 0 & 0 & -1 & 0 \\
                   0 & -1 & 0 & 1 & 0 & 0 & 0 & 0 & 0 & 0 & 0 & 0 & 0 & 0 & 0 & 0 & 0 & 1 \\
                   0 & -1 & 1 & 0 & 0 & 0 & 0 & 0 & 0 & 0 & 0 & 0 & 0 & 0 & 0 & 0 & -1 & -1 \\
                   0 & 0 & -1 & 1 & 0 & 0 & 0 & 0 & 0 & 0 & 0 & 0 & 0 & 0 & 0 & 0 & 0 & -1 \\
                   0 & 0 & -1 & 0 & 1 & 0 & 0 & 0 & 0 & 0 & 0 & 0 & 0 & 0 & 0 & 0 & 0 & 0 \\
                   0 & 0 & 0 & -1 & 0 & 1 & 0 & 0 & 0 & 0 & 0 & 0 & 0 & 0 & 0 & 0 & 0 & 0 \\
                   0 & 0 & 0 & -1 & 0 & 0 & 1 & 0 & 0 & 0 & 0 & 0 & 0 & 0 & 0 & 0 & 0 & 0 \\
                   0 & 0 & 0 & 0 & -1 & 1 & 0 & 0 & 0 & 0 & 0 & 0 & 0 & 0 & 0 & 0 & 0 & 0 \\ \hline
                   0 & 0 & 0 & 0 & 0 & 0 & 0 & 1 & -1 & 0 & -1 & 0 & 0 & 0 & 0 & 0 & 0 & 0 \\
                   0 & 0 & 0 & 0 & 0 & 0 & 0 & 0 & 0 & 1 & -1 & -1 & 0 & 0 & 0 & 0 & 0 & 0 \\
                  \end{array}\right]$$
\end{tiny}
The gradient $d_0$ and curl $d_1$ are 
\begin{tiny}
$$ d_0 = \left[
                  \begin{array}{ccccccc}
                   -1 & 0 & 1 & 0 & 0 & 0 & 0 \\
                   -1 & 1 & 0 & 0 & 0 & 0 & 0 \\
                   0 & -1 & 0 & 1 & 0 & 0 & 0 \\
                   0 & -1 & 1 & 0 & 0 & 0 & 0 \\
                   0 & 0 & -1 & 1 & 0 & 0 & 0 \\
                   0 & 0 & -1 & 0 & 1 & 0 & 0 \\
                   0 & 0 & 0 & -1 & 0 & 1 & 0 \\
                   0 & 0 & 0 & -1 & 0 & 0 & 1 \\
                   0 & 0 & 0 & 0 & -1 & 1 & 0 \\
                  \end{array}
                  \right], 
  d_1 = \left[ 
                  \begin{array}{ccccccccc}
                   1 & -1 & 0 & -1 & 0 & 0 & 0 & 0 & 0 \\
                   0 & 0 & 1 & -1 & -1 & 0 & 0 & 0 & 0 \\
                  \end{array}
                  \right]  \; . $$
\end{tiny}
The characteristic polynomial of $D$ is $p_D(x) = x^{18}-24 x^{16}+242 x^{14}-1334 x^{12}+4377 x^{10}-8706 x^8+10187 x^6-6370 x^4+1624 x^2$. 
The positive eigenvalues of $D$ are $0.92, 1.05, 1.41, 1.69, 1.78, 2.00, 2.15, 2.38$.
The kernel of $L_0$ is spanned by $[1,1,1,1,1,1,1]^T$, the kernel of $L_1$ is spanned by
$[1,-1,-3,2,-5,8,-8,0,8]^T$. The Laplace-Beltrami operator $L$ of $G$ has three blocks:
\begin{tiny}
$$  L = \left[
                  \begin{array}{ccccccc|ccccccccc|cc}
                   2  & -1 & -1 &  0 &  0 &  0 &  0 & 0 & 0 & 0 & 0 & 0 & 0 & 0 & 0 & 0 & 0 & 0 \\
                   -1 & 3  & -1 & -1 &  0 &  0 &  0 & 0 & 0 & 0 & 0 & 0 & 0 & 0 & 0 & 0 & 0 & 0 \\
                   -1 & -1 & 4 &  -1 & -1 &  0 &  0 & 0 & 0 & 0 & 0 & 0 & 0 & 0 & 0 & 0 & 0 & 0 \\
                   0  & -1 & -1 &  4 &  0 & -1 & -1 & 0 & 0 & 0 & 0 & 0 & 0 & 0 & 0 & 0 & 0 & 0 \\
                   0  & 0  & -1 &  0 &  2 & -1 &  0 & 0 & 0 & 0 & 0 & 0 & 0 & 0 & 0 & 0 & 0 & 0 \\
                   0  & 0  &  0 & -1 & -1 &  2 &  0 & 0 & 0 & 0 & 0 & 0 & 0 & 0 & 0 & 0 & 0 & 0 \\
                   0  & 0  &  0 & -1 &  0 &  0 &  1 & 0 & 0 & 0 & 0 & 0 & 0 & 0 & 0 & 0 & 0 & 0 \\ \hline
                   0  & 0  &  0 &  0 &  0 &  0 &  0 & 3 & 0 & 0 & 0 & -1 & -1 & 0 & 0 & 0 & 0 & 0 \\
                   0  & 0  &  0 &  0 &  0 &  0 &  0 & 0 & 3 & -1 & 0 & 0 & 0 & 0 & 0 & 0 & 0 & 0 \\
                   0  & 0  &  0 &  0 &  0 &  0 &  0 & 0 & -1 & 3 & 0 & 0 & 0 & -1 & -1 & 0 & 0 & 0 \\
                   0  & 0  &  0 &  0 &  0 &  0 &  0 & 0 & 0 & 0 & 4 & 0 & -1 & 0 & 0 & 0 & 0 & 0 \\
                   0  & 0  &  0 &  0 &  0 &  0 &  0 & -1 & 0 & 0 & 0 & 3 & 1 & -1 & -1 & 0 & 0 & 0 \\
                   0  & 0  &  0 &  0 &  0 &  0 &  0 & -1 & 0 & 0 & -1 & 1 & 2 & 0 & 0 & -1 & 0 & 0 \\
                   0  & 0  &  0 &  0 &  0 &  0 &  0 & 0 & 0 & -1 & 0 & -1 & 0 & 2 & 1 & 1 & 0 & 0 \\
                   0  & 0  &  0 &  0 &  0 &  0 &  0 & 0 & 0 & -1 & 0 & -1 & 0 & 1 & 2 & 0 & 0 & 0 \\
                   0  & 0  &  0 &  0 &  0 &  0 &  0 & 0 & 0 & 0 & 0 & 0 & -1 & 1 & 0 & 2 & 0 & 0 \\ \hline
                   0  & 0  &  0 &  0 &  0 &  0 &  0 & 0 & 0 & 0 & 0 & 0 & 0 & 0 & 0 & 0 & 3 & 1 \\
                   0  & 0  &  0 &  0 &  0 &  0 &  0 & 0 & 0 & 0 & 0 & 0 & 0 & 0 & 0 & 0 & 1 & 3 \\
                  \end{array}
                  \right] \; . $$
\end{tiny}
The last block $L_2 = \left[ \begin{array}{cc} 3 & 1 \\ 1 & 3 \end{array} \right]$ acts on functions on
triangles. That the diagonal entries are $3$ follows from ${\rm deg}_p(x) = L_p(x,x)-(p+1)$ for $p=2$
and the fact that there are no tetrahedra so that the degree is zero for every triangle. 

\section{Mathematica code}
We illustrate how brief the procedure building the Dirac operator from a graph can be. 
The source TeX File can be accessed, when opening the source of this file on the ArXiv. 
Source code to earlier papers can be found on my website. 
\pagebreak

\begin{tiny}
\lstset{language=Mathematica} \lstset{frameround=fttt}
\begin{lstlisting}[frame=single]
Cliques[K_,k_]:=Module[{n,u,m,s,V=VertexList[K],W=EdgeList[K],U,r={}},
 s=Subsets[V,{k,k}]; n=Length[V]; m=Length[W]; 
 Y=Table[{W[[j,1]],W[[j,2]]},{j,Length[W]}]; If[k==1,r=V,If[k==2,r=Y,
 Do[u=Subgraph[K,s[[j]]]; If[Length[EdgeList[u]]==Binomial[k,2],
 r=Append[r,VertexList[u]]],{j,Length[s]}]]];r];
Dirac[s_]:=Module[{a,b,q,l,n,v,m,R,t,d},q=VertexList[s];n=Length[q];
 d=Table[{{0}},{p,n-1}]; l=Table[{},{p,n}]; v=Table[0,{p,n}]; m=n; 
 Do[If[m==n,l[[p]]=Cliques[s,p];v[[p]]=Length[l[[p]]]; 
 If[v[[p]]==0,m=p-2]],{p,n}]; t=Sum[v[[p]],{p,n}]; 
 b=Prepend[Table[Sum[v[[p]],{p,1,k}],{k,Min[n,m+1]}],0]; 
 R=Table[0,{t},{t}]; If[m>0, d[[1]] = Table[0,{j,v[[2]]},{i,v[[1]]}]; 
   Do[d[[1,j,l[[2,j,1]]]]=-1,{j,v[[2]]}]; 
   Do[d[[1,j,l[[2,j,2]]]]=1,{j,v[[2]]}]]; 
 Do[ If[m>=p,d[[p]]=Table[0,{j,v[[p+1]]},{i,v[[p]]}]; 
   Do[a=l[[p+1,i]];Do[d[[p,i,Position[l[[p]],
   Delete[a,j]][[1,1]]]]=(-1)^j,{j,p+1}],{i,v[[p+1]]}]],{p,2,n-1}];
 Do[If[m>=p,Do[R[[b[[p+1]]+j,b[[p]]+i]]=d[[p,j,i]],
   {i,v[[p]]},{j,v[[p+1]]}]],{p,n-1}]; 
 R+Transpose[R]]; 
s={1->2,2->3,3->1,3->4,4->2,3->5,5->6,6->4,4->7};
s=UndirectedGraph[Graph[s]]; DD=Dirac[s] 
\end{lstlisting}
\end{tiny}


\begin{thebibliography}{10}

\bibitem{josellisknill}
F.~Josellis and O.~Knill.
\newblock A {L}usternik-{S}chnirelmann theorem for graphs.
\newblock \\http://arxiv.org/abs/1211.0750, 2012.

\bibitem{randomgraph}
O.~Knill.
\newblock The dimension and {Euler} characteristic of random graphs.
\newblock {\\}http://arxiv.org/abs/1112.5749, 2011.

\bibitem{cherngaussbonnet}
O.~Knill.
\newblock A graph theoretical {Gauss-Bonnet-Chern} theorem.
\newblock {\\}http://arxiv.org/abs/1111.5395, 2011.

\bibitem{elemente11}
O.~Knill.
\newblock A discrete {Gauss-Bonnet} type theorem.
\newblock {\em Elemente der Mathematik}, 67:1--17, 2012.

\bibitem{poincarehopf}
O.~Knill.
\newblock A graph theoretical {Poincar\'e-Hopf} theorem.
\newblock {\\} http://arxiv.org/abs/1201.1162, 2012.

\bibitem{indexformula}
O.~Knill.
\newblock An index formula for simple graphs \hfill.
\newblock {\\}http://arxiv.org/abs/1205.0306, 2012.

\bibitem{indexexpectation}
O.~Knill.
\newblock On index expectation and curvature for networks.
\newblock {\\}http://arxiv.org/abs/1202.4514, 2012.

\bibitem{knillmckeansinger}
O.~Knill.
\newblock {The McKean-Singer Formula in Graph Theory}.
\newblock {\\}http://arxiv.org/abs/1301.1408, 2012.

\bibitem{knillcalculus}
O.~Knill.
\newblock {The theorems of Green-Stokes,Gauss-Bonnet and Poincare-Hopf in Graph
  Theory}.
\newblock {\\}http://arxiv.org/abs/1201.6049, 2012.

\bibitem{brouwergraph}
O.~Knill.
\newblock A {Brouwer} fixed point theorem for graph endomorphisms.
\newblock {\em Fixed Point Theory and Applications}, 85, 2013.

\bibitem{cauchybinet}
O.~Knill.
\newblock Cauchy-{B}inet for pseudo determinants.
\newblock {\\}http://arxiv.org/abs/1306.0062, 2013.

\bibitem{diracannouncement}
O.~Knill.
\newblock An integrable evolution equation in geometry.
\newblock {\\} http://arxiv.org/abs/1306.0060, 2013.

\end{thebibliography}

\end{document}